 \UseRawInputEncoding 
\RequirePackage{fix-cm}

\documentclass[smallextended]{svjour3}       
\smartqed  
\usepackage[square, comma, sort&compress, numbers]{natbib}
\usepackage{graphicx}
\usepackage[misc]{ifsym}
\usepackage{amssymb}
\usepackage{mathrsfs}
\usepackage{bm}
\usepackage{float}  
\usepackage{subfigure}  

\usepackage{amsmath}
\allowdisplaybreaks[4]
\usepackage{array}
\usepackage{graphicx}
\usepackage{cases}
\usepackage{color}
\usepackage{multirow}
\usepackage{picinpar}
\graphicspath{ {./images/} }

\usepackage[all]{xy}
\usepackage[colorlinks,
            linkcolor=blue,
            anchorcolor=green,
            citecolor=blue
            ]{hyperref}

\def\N{\mathcal{N}}

\def\N{{\mathbb N}}

\def\R{{\mathbb R}}

\def\bfe{{\bf 1}}
\def\bfe{{\bf 1}}

\def\bfe{{\bf 1}}

\def\bfv{{\bf v}}
\def\bfx{{\bf x}}
\def\bfw{{\bf w}}
\def\bfy{{\bf y}}

\def\bfs{{\bf s}}
\def\bfu{{\bf u}}

\def\R{\mathbb{R}}

\newcommand{\ba}{\begin{array}}
\newcommand{\ea}{\end{array}}
\newcommand{\be}{\begin{equation}}
\newcommand{\ee}{\end{equation}}

\newcommand{\RNum}[1]{\lowercase\expandafter{\romannumeral #1\relax}}

\begin{document}
\title{Proximal Operator and Optimality Conditions for Ramp Loss SVM
}

\titlerunning{POOC for $L_r$-SVM }        
\author{Huajun Wang$^1$ \and
        Yuanhai Shao$^2$  \and
        Naihua Xiu$^1$
}

\institute{Huajun Wang \at huajunwang@bjtu.edu.cn\and
           Yuanhai Shao \at shaoyuanhai@hainanu.edu.cn\and
           Naihua Xiu \at nhxiu@bjtu.edu.cn\and
           1. Department of Mathematics, School of Science,
Beijing Jiaotong University, Beijing \at 100044, People's Republic of China\and
 2. School of Management, Hainan University, Haikou, \at 570000, People's Republic of China
           }

\date{Received: date / Accepted: date}

\maketitle

\begin{abstract}

Support vector machines with ramp loss (dubbed as $L_r$-SVM) have attracted wide attention due to the boundedness of ramp loss. However, the corresponding optimization problem is non-convex and the given Karush-Kuhn-Tucker (KKT) conditions are only the necessary  conditions. To enrich the optimality theory of $L_r$-SVM and go deep into its statistical nature, we first introduce and analyze the proximal operator for ramp loss, and then establish a stronger optimality conditions: P-stationarity, which is proved to be the first-order necessary and sufficient conditions for local minimizer of $L_r$-SVM. Finally, we define the $L_r$ support vectors based on the concept of P-stationary point, and show that all $L_r$ support vectors fall into the support hyperplanes, which possesses the same feature as the one of hard margin SVM.

\end{abstract}
\vskip 2mm
\section{Introduction}
 \vskip 3mm

Support vector machines (SVM) were first introduced by Vapnik and Cortes \cite{CV95} and have been  widely
 applied into many fields, including text and image classification \cite{ZY08,LZ20}, disease detection \cite{GG19,PH07}, etc.
The decision hyperplane of SVM classifier, $\langle  \bfw , \bfx\rangle+b=0$ with $ \bfw \in \R^n$ and $b\in \R$, is trained from data set $\{(\bfx_i,y_i), i\in\N_m\}$ where $\bfx_i\in \R^n$, $y_i\in \{-1,1\}$ and $\N_m:=\{1,2,\cdots,m\}$ by optimizing
the following problem
\begin{eqnarray}\label{SM-SVM}
\min_{ \bfw \in \R^n,b\in \R}~~\frac{1}{2}\Vert \bfw \Vert^{2}+C\sum_{i=1}^{m}\ell_h(1-y_{i}( \langle \bfw, \bfx_i\rangle +b)),
\end{eqnarray}
where $C>0$ is a penalty parameter and $\ell_h(t)=\max\{t,0\}$ is the hinge  loss, which has no cost for $t<0$, but pays linear cost for $t\geq0$. The cost given to the outliers by the hinge loss is quite huge since it is unbounded function \cite{HS14,DE15,MQ20,WF17,KS09,LS05}.  To remove the impact of outliers, one method for increasing the robustness of SVM is to use the ramp loss \cite{ST03} (see Fig. 1 (a)).  The ramp loss is defined as follows,
\begin{eqnarray}\label{Rloss}
\ell_{r}(t)=
\begin{cases}
1,& t\geq1,\\
t,& 0\leq t<1,\\
0,& t<0,
\end{cases}
\end{eqnarray}
which is also called truncated hinge loss in \cite{YY07}. It has no cost for $t<0$, but it pays linear cost for $0\leq t< 1$ and a fixed cost at 1 for $t\geq1$. Authors in \cite{RF06a,RF06b} also extended and studied ramp loss by replacing all constant 1 in (\ref{Rloss}) by adjustable parameter $\mu$ $(\mu>0)$.  Since the theoretical results have no essential difference between the fixed 1 and parameter $\mu$ for ramp loss, for simplicity of the proof, we use formula \eqref{Rloss} directly. Bartlett et al. \cite{BM02} investigated some of the theoretical
properties of the ramp loss, while SVM with ramp loss  (dubbed as $L_r$-SVM) was first proposed by Shen et al. \cite{ST03}. Replacing $\ell_{h}(\cdot)$ by $\ell_{r}(\cdot)$ in (\ref{SM-SVM}), the optimization problem of $L_r$-SVM is 
\begin{eqnarray} \label{R-SVM}
\min_{ \bfw \in {\R}^n, b\in {\R}} f_r( \bfw ;b):=\frac{1}{2}\Vert \bfw  \Vert^2+CL_r({\bfe}-A \bfw -b{\bf{y}}),
\end{eqnarray}
where $A:=[y_1\bfx_{1}~y_2\bfx_{2}~\cdots~y_m\bfx_{m}]^\top\in {\R}^{m\times n}, \bfy:=(y_{1},y_2,\cdots,y_{m})^{\top}\in {\R}^{m}, {\bfe} :=(1,1,\cdots,1)^{\top}\in {\R}^{m}$, $L_r(\bfu):=\sum_{i=1}^m \ell_{r}(u_i)$ with $\bfu:=(u_1,u_2,\cdots,u_m)^\top={\bfe}-A \bfw -b{\bf{y}}\in\R^m$, which computes the sum of positive elements for all $u_i<1$ and the number of elements for all $u_i\geq 1$ in $\bfu$. Due to the non-convexity of ramp loss, the $L_r$-SVM (\ref{R-SVM}) is a non-differentiable non-convex optimization problem. Thus, to find a satisfying solution, scholars have paid many efforts in optimality conditions and algorithms for $L_r$-SVM.

 In recent decades, on the one hand, some approaches such as convex relaxation and equivalent reformulation as mixed integer nonlinear programming (MINP) problem were proposed for dealing with $L_r$-SVM. Xu et al. \cite{Xu06} reformulated $L_r$-SVM as a semidefinite programming problem by convex relaxation which was solved by the MATLAB software package SDPT3.  Brooks \cite{BJ11}  transformed  $L_r$-SVM into a MINP problem, and then proposed a branch and bound algorithm to solve it. Carrizos et al. \cite{CN14} developed  heuristic method to handle the  MINP problem of $L_r$-SVM on large datasets.

On the other hand, Collobert et al. \cite{RF06a} translated $L_r$-SVM (\ref{R-SVM}) into an equivalent  difference of convex functions (DC) programming and took the subdifferential of $f_r(\bfw;b)$ to establish the first-order necessary conditions: Karush-Kuhn-Tucker (KKT) conditions, and used the concave-convex
procedure (CCCP) \cite{YR03} for solving the DC programming.  To improve the computational speed, Wang et al. \cite{ZS09} designed an efficient working set selection strategy based on KKT conditions and then used the CCCP with working set for solving the DC programming.
 Since the KKT conditions provided effective characterization of optimal solution to $L_r$-SVM (\ref{R-SVM}), this class of approaches continued to be widely studied in theory as well as algorithms. For more details, see, e.g., \cite{EB11,WL20,LW19,SS18} and references therein.

 A natural question arises. Is there stronger first-order optimality conditions of $L_r$-SVM (\ref{R-SVM}) to provide more effective characterization of optimal solution?
It is this question that motivates the work in our paper.
 The main results of this paper are summarized as follow: (i) We derive the explicit expression of proximal operator of ramp  loss.  (ii) We introduce a novel optimality conditions: P-stationarity, which is proved to be the necessary and sufficient conditions for local minimizer of $L_r$-SVM. (iii) We prove that the set of the support vectors defined by P-stationary point fall into the support hyperplanes.


This paper is organized as follows. In Section 2, we derive the explicit expression of proximal operator of ramp  loss. In Section 3, we introduce a concept of P-stationary point, and reveal the relationship between P-stationary point and local/global minimizer of $L_r$-SVM. In Section 4, we introduce $L_r$ support vectors based on the P-stationary point and discuss its properties. Conclusions are made in Section 5.
\vskip 1mm
\section{Proximal operator for ramp loss}
In this section,  we derive the explicit expression of proximal operator of ramp loss, which will be used to study the new first-order optimality conditions of $L_r$-SVM in the next section.
{\bf \subsection{ $\ell_{r}$  proximal operator }}
We first give the definition of $\ell_{r}$ proximal operator in one-dimensional case.
\begin{definition}
($\ell_{r}$ proximal operator) For any given $\gamma, C>0$ and ${s}\in \R$, the proximal operator of $\ell_{r}(v)$ (dubbed as $\ell_{r}$ proximal operator) is defined as
\begin{eqnarray} \label{proxmalL}
\text{prox}_{\gamma C \ell_{r}}({ s})=\arg \min_{v\in {\R}}~ C
\ell_{r}(v)+\frac{1}{2\gamma}(v-s)^2.
\end{eqnarray}\end{definition}

The following two propositions state that the $\ell_{r}$ proximal operator admits a closed form solution for  $0<\gamma C<2$ or $\gamma C\geq2$.
\begin{proposition}\label{pro1} {\rm(Solution to $\ell_{r}$ proximal operator for $0<\gamma C<2$) \label{proximalLma1} For any given $\gamma, C>0$ and $0<\gamma C<2$, the solution to $\ell_{r}$ proximal operator  at $s \in \R$ is given as
}\end{proposition}
\be \label{prox1}
\text{prox}_{\gamma C\ell_{r}}(s):=
\begin{cases}
s,& s>1+\frac{\gamma C}{2},\\
s ~\text{or}~ s-\gamma C,& s=1+\frac{\gamma C}{2},\\
{s-\gamma C} ,&\gamma C\leq s<1+\frac{\gamma C}{2},\\
{0} ,&0< s< \gamma C,\\
s,&  \ s\leq0.
\end{cases}
\ee

\begin{proposition}{ \label{prox2}\rm(Solution to $\ell_{r}$ proximal operator for $\gamma C\geq2$) \label{proximalLma1} For any given $\gamma, C>0$ and $\gamma C\geq2$, the solution to $\ell_{r}$ proximal operator  at $s \in \R$ is given as
}\end{proposition}
\be \label{prox12}
\text{prox}_{\gamma C\ell_{r}}(s):=
\begin{cases}
s,& s>\sqrt{2\gamma C},\\
 s~\text{or}~0,& s=\sqrt{2\gamma C},\\
0 ,&0< s<\sqrt{2\gamma C},\\
s,&  \ s\leq0.
\end{cases}
\ee

{\bf\subsection{ $L_{r}$ proximal operator }}
Based on separate property of $L_{r}(\cdot)$, we extend \eqref{prox1} and \eqref{prox12} to  multi-dimensional case.
\vskip 2mm
\begin{definition}($L_{r}$ proximal operator) \label{peox_define}For any given $\gamma, C>0$, the proximal operator of $L_{r}(\bfv)$ (dubbed as $L_{r}$ proximal operator) at ${\bf s}=(s_1,s_2,\cdots,s_m)^{\top}\in\R^m$ is defined as
\begin{eqnarray} \label{proxmalLL}
\text{prox}_{\gamma CL_{r}}({\bf s})=\arg\min_{\bfv\in {\R}^m}~ C
L_{r}(\bfv)+\frac{1}{2 \gamma}\|\bfv-{\bf s}\|^2.
\end{eqnarray}
\end{definition}

The following proposition states that the $L_{r}$ proximal operator admits a closed form solution for two cases: $0<\gamma C<2$ or $\gamma C\geq2$.
\begin{proposition}\label{pro3}{\rm(Solution to $L_{r}$ proximal operator) \label{multi}For a given $\gamma, C>0$, the solution to $L_{r}$ proximal operator at ${\bf s}=(s_1,s_2,\cdots,s_m)^{\top}\in\R^m$ is given as
\begin{eqnarray}\label{Prox_L}
\text{prox}_{\gamma CL_{r}}({ \bfs}):=
\left[\begin{array}{c}
\text{prox}_{\gamma C\ell_{r}}(s_1)\\
\vdots\\
{ \text{prox}_{\gamma C\ell_{r}}(s_m)}
\end{array}\right],\end{eqnarray}
where $\text{prox}_{\gamma C\ell_{r}}(s_i)$ takes formula (\ref{prox1}) as $0<\gamma C<2$ or (\ref{prox12}) as $\gamma C\geq2$.
}\end{proposition}

\noindent {\bf \emph{Proof}} It follows from (\ref{proxmalLL}) that $[\text{prox}_{\gamma CL_{r}}({ \bfs})]_i=\text{prox}_{\gamma C\ell_{r}}({ s_i})$, where
\begin{align*}
\text{prox}_{\gamma C\ell_{r}}(s_i)=\text{arg}\min_{v\in {\R}} C\ell_{r}(v)+ \frac{1}{2\gamma}(v-s_i)^{2}, i\in \N_m.
\end{align*}
 Using (\ref{prox1}) or (\ref{prox12}) completes the proof.  \hfill  $ \Box $\vskip 2mm

\section{ First-order optimality conditions}

In this section, we develop a first-order necessary and sufficient optimality conditions for (\ref{R-SVM}). To proceed this, we introduce a variable ${\bfu}\in {\R}^{m}$ to equivalently  reformulate (\ref{R-SVM}) as
\begin{eqnarray} \label{RR-SVM}
&\underset{ \bfw \in {\R}^n,b \in \R,{\bfu} \in {\R}^m}{\min} & \frac{1}{2} \Vert \bfw  \Vert^2+CL_r(\bfu)\\ \nonumber
&\mbox{s.t.} &{\bfu}+A  \bfw +b\bfy={\bfe}.
\end{eqnarray}
Now let us define some notation
\begin{eqnarray}\label{BH}
B:=[A~\bfy]\in\R^{m\times(n+1)},~~~~
H:=\begin{bmatrix}
I_{n\times n} & {\bf0}  \\
{\bf0} & 0
\end{bmatrix}B^\dag,
\end{eqnarray}
 where  $B^\dag\in\R^{(n+1)\times m}$ is the generalized inverse of $B$ and $I_{n\times n}$ is the identity matrix with order $n$. Denote
$\lambda_H:= \lambda_{\max}(H^\top H)$ where $\lambda_{\max}(H^\top H)$ is the maximum eigenvalue of $H^\top H.$
\begin{definition}
(P-stationary point of (\ref{RR-SVM})) For a  given $C>0$, we call $({ \bfw^* }; b^*;  \bfu^* )$ is a proximal stationary (P-stationary) point of (\ref{RR-SVM}) if  there exists a Lagrangian multiplier vector $ {\bm{\lambda}^*}\in \R^m$ and a constant $\gamma>0$ such that
\be\label{PH-sta}
\left\{
\begin{array}{rll}   \bfw^* + A^{\top}{\bm{\lambda}^*} &=& {\bf 0},
\\ \langle \bfy, {\bm{\lambda}^*} \rangle&=&  0,
\\  \bfu^*+A   \bfw^*+b^*\bfy&=&
{\bfe},
\\ \text{prox}_{\gamma C L_{r}}(\bfu^*-\gamma{\bm{\lambda}^*})&\ni&\bfu^*.
\end{array}\right.
\ee
\end{definition}

Based on the above definition, we obtain the desired result in this section.

\begin{theorem}{\rm\label{theom1} (First-order necessary optimality conditions). Let $B$ be full  column rank. For a given $C>0$, if $( \bfw ^*; b^*; {\bfu}^*)$ is a global minimizer of (\ref{RR-SVM}), then it is a P-stationary point with $0<\gamma< 1/\lambda_H$.

 }\end{theorem}

 \begin{theorem}{\rm\label{theom2} (First-order sufficient optimality conditions). For a given $C>0$, if $( \bfw ^*; b^*; {\bfu}^*)$ with ${\bm{\lambda}}^*\in {\R}^m$ and $\gamma>0$ is a P-stationary point, then it is a local minimizer of (\ref{RR-SVM}).
 }\end{theorem}

\vskip 2mm
At the end of this section, we analyze the  following relationship between the P-stationary point  and the KKT point given by \cite{RF06a}. Based on the problem of \eqref{RR-SVM} and the necessary conditions results in \cite{RF06a}, if $( {\bf\overline{w}} ; \overline{b}; {\bf\overline{u}} )$ is a local minimizer of (\ref{RR-SVM}), then we have the following KKT conditions,

\be\label{PH-stakkt}
\left\{
\begin{array}{rll}   \overline{\bfw} + A^{\top}{\bm{\overline{\lambda}}} &=& {\bf 0},
\\ \langle \bfy, {\bm{\overline{\lambda}}} \rangle&=&  0,
\\  \overline{\bfu}+A   \overline{\bfw}+\overline{b}\bfy&=&
{\bfe},
\\ C\partial L_{r}(\overline{\bfu})+\bm{\overline{\lambda}}&\ni&{\bf 0},
\end{array}\right.
\ee
where ${\bm{\overline{\lambda}}} \in {\R}^m$ is a multiplier vector, $( {\bf\overline{w}} ; \overline{b}; {\bf\overline{u}} )$ is called the KKT point of \eqref{RR-SVM}, and $\partial L_r({\bf\overline{u}})=(\partial \ell_r(\overline{u}_1),\partial \ell_r(\overline{u}_2),\cdots,\partial \ell_r(\overline{u}_m))^\top\in \R^m$ is the  subdifferential of $L_r({\bf\overline{u}})$ with 
\be\nonumber
\partial\ell_r(\overline{u}_i)
\begin{cases}
\in[0,1],& \overline{u}_i\in\{0,1\},\\
=1,& \overline{u}_i\in(0,1),\\
=0,&\overline{u}_i<0~\text{or}~\overline{u}_i>1 ,
\end{cases}~~i\in \N_m.
\ee
Furthermore, from the above formula and last formula of \eqref{PH-stakkt}, we have
\begin{equation}\label{KKT-lam}
\textbf{0}\in C\partial
L_r({\bf\overline{u}})+{\bm{\overline{\lambda}}} \Leftrightarrow \left\{\bm{\overline{\lambda}}\in {\R}^m:
\overline{\lambda}_i
\begin{cases}
 \in[-C,0], &\overline{u}_i\in\{0,1\},\\
=-C,&\overline{u}_i\in(0,1),
\\ =0,&\overline{u}_i<0~\text{or}~\overline{u}_i>1.
 \end{cases}
 \right.
\end{equation}


\begin{theorem}{\rm\label{P-KKT}  For a given $C>0$, if $( \bfw ^*; b^*; {\bfu}^* )$ with $\bm{{\lambda}}^*\in\R^m$ and $\gamma>0$ is a P-stationary point of (\ref{RR-SVM}), then it is also a KKT point of (\ref{RR-SVM}),  but the converse does not hold.
}\end{theorem}
\noindent {\bf \emph{Proof}} The former conclusion  is obvious since the P-stationary point is  the local minimizer of (\ref{RR-SVM}). The latter conclusion is from the following counterexample. Consider the training set with the positive vectors $\bfx_1=(3,3)^\top$, $\bfx_2=(6,-2)^\top$ and negative
vector $\bfx_3=(1,1)^\top$. Namely,
\begin{eqnarray*}
A:=[y_1\bfx_{1}~y_2\bfx_{2}~y_3\bfx_{3}]^\top=
\begin{bmatrix}
3&6&-1\\
3&-2&-1
\end{bmatrix}^\top,~~\bfy=\begin{bmatrix}
 1,&1,&-1
\end{bmatrix}^\top.
\end{eqnarray*}
For a given $C=0.25$, we can verify that
\begin{eqnarray*}
{\bf\overline{w}}&=&(0.5,0.5)^\top,~~~~~~~~~~~~~~~ \overline{b}=-2,~~~~~~~~~ {\bf\overline{u}}=(0,1,0)^\top,\\
\bm{\overline{\lambda}}&=&(-0.25,0,-0.25)^\top, ~~~~\partial
L_r({\bf\overline{u}})=(1,0,1)^\top
\end{eqnarray*}
 satisfy \eqref{PH-stakkt}. This means  $( {\bf\overline{w}} ; \overline{b}; {\bf\overline{u}} )$ with $\bm{\overline{\lambda}}$ is a KKT point of \eqref{RR-SVM}.  Particularly, for $i=2$, we have $\overline{u}_2=1$ and $\overline{\lambda}_2=0$. However, for $0<\gamma<8$ and $C=0.25$, i.e., $0<\gamma C<2$, from \eqref{prox1} and $\overline{u}_2-\gamma\overline{\lambda}_2\in (0, 1+\frac{\gamma C}{2})$, we obtain $\text{prox}_{\gamma C\ell_{r}}(\overline{u}_2-\gamma\overline{\lambda}_2)=0~ \text{or}~ \overline{u}_2-\gamma C$ but the both does not equal $\overline{u}_2$. For $\gamma\geq8$ and $C=0.25$, i.e., $\gamma C\geq2$, from \eqref{prox12} and $\overline{u}_2-\gamma\overline{\lambda}_2\in (0, \sqrt{2\gamma C})$, we get $\text{prox}_{\gamma C\ell_{r}}(\overline{u}_2-\gamma\overline{\lambda}_2)= 0\neq \overline{u}_2$. To sum up, we have $\overline{u}_2\notin\text{prox}_{\gamma C\ell_{r}}(\overline{u}_2-\gamma\overline{\lambda}_2)$, which means $( {\bf\overline{w}} ; \overline{b}; {\bf\overline{u}} )$  is not a P-stationary point of (\ref{RR-SVM}). \hfill $\Box $

\section{$L_r$ support vectors}\label{sec:l01-sv}

In this section, we define the support vectors  for $L_r$-SVM by P-stationary point of (\ref{RR-SVM}), which are named as $L_r$ support vectors since they are selected by the $L_r$ proximal operator. Before proceeding, we first review  the definitions of support vectors of hard margin SVM, hinge loss SVM and ramp loss SVM, respectively.

 {\bf i) Support vectors of hard margin SVM\cite{CV95}:} Let (${\bf\widetilde{w}}; \widetilde{b}$) be a global minimizer of hard margin SVM and  $\bm{\widetilde{\alpha}}=(\widetilde{\alpha}_1,\widetilde{\alpha}_2,\cdots,\widetilde{\alpha}_m)^\top$ with $\widetilde{\alpha}_i\geq0$  be a solution of its dual problem. Then the  global minimizer ${\bf\widetilde{w}}$ satisfies
 \begin{eqnarray}\nonumber
\label{ww1}{\bf\widetilde{w}} =\underset{i \in \widetilde{J}^{*}}{\sum}\widetilde{\alpha}_i y_i\bfx_{i},
\end{eqnarray}
where $\widetilde{J}^{*}:= \{i\in \N_m:\widetilde{\alpha}_i>0\}$. The training vectors $\{\bfx_i, i\in \widetilde{J}^*\}$ are called support vectors. For any $i\in \widetilde{J}^{*}$, we have
 \begin{eqnarray}\nonumber
 y_i(\langle{\bf\widetilde{w}},\bfx_{i}\rangle+\widetilde{b})=1.
\end{eqnarray}

 {\bf ii) Support vectors of hinge loss SVM\cite{CV95}:} Let (${\bf\widehat{w}}; \widehat{b}$) be a  global minimizer of hinge loss SVM and  $\bm{\widehat{\alpha}}=(\widehat{\alpha}_1,\widehat{\alpha}_2,\cdots,\widehat{\alpha}_m)^\top$ with $\widehat{\alpha}_i\in[0, C]$ be a solution of its dual problem. Then the  global minimizer $\bf\widehat{w}$ satisfies
 \begin{eqnarray}\nonumber
\label{ww1}{\bf\widehat{w}} =\underset{i \in \widehat{J}^{*}}{\sum} \widehat{\alpha}_iy_i\bfx_{i},
\end{eqnarray}
where $\widehat{J}^{*}:=\{i\in \N_m:\widehat{\alpha}_i\in(0,C]\}$. The training vectors $\{\bfx_i, i\in \widehat{J}^*\}$ are called support vectors. For any $i\in\widehat{J}^{*}$, we have
 \begin{eqnarray}\nonumber
\begin{cases} y_i(\langle{\bf\widehat{w}},\bfx_{i}\rangle+\widehat{b})=1, ~ i \in \{i\in \N_m:\widehat{\alpha}_i\in(0,C)\},\\
y_i(\langle{\bf\widehat{w}},\bfx_{i}\rangle+\widehat{b})\leq1, ~ i \in \{i\in \N_m:\widehat{\alpha}_i=C\}.
\end{cases}
\end{eqnarray}

{\bf iii) Support vectors of ramp loss SVM\cite{RF06a}:} Let (${\bf\overline{w}}; \overline{b}; \overline{\bfu}$) and $\bm{\overline{\lambda}}\in \R^m$ with $\overline{\lambda}_i\in[-C,0]$ satisfy  KKT conditions. From \eqref{KKT-lam} and the first equation of \eqref{PH-stakkt}, the KKT point ${\bf\overline{w}}$ satisfies
 \begin{eqnarray}\nonumber
\label{ww1}{\bf\overline{w}} =-\underset{i \in \overline{J}^{*}}{\sum} \overline{\lambda}_iy_i\bfx_{i},
\end{eqnarray}
where $\overline{J}^{*}:=\{i\in \N_m:\overline{\lambda}_i\in[-C,0)$. The training vectors $\{\bfx_i, i\in \overline{J}^*\}$ are called support vectors. For any $i\in \overline{J}^{*}$, from \eqref{KKT-lam} and the third equation of \eqref{PH-stakkt}, we have
 \begin{eqnarray}\nonumber
\begin{cases} y_i(\langle{\bf\overline{w}},\bfx_{i}\rangle+\overline{b})\in\{0,1\}, ~ i \in \{i\in \N_m:\overline{\lambda}_i\in(-C,0)\},\\
y_i(\langle{\bf\overline{w}},\bfx_{i}\rangle+\overline{b})\in[0,1], ~ i \in \{i\in \N_m:\overline{\lambda}_i=-C\}.
\end{cases}
\end{eqnarray}

The above results show that the hard margin SVM and hinge loss SVM  define the support vectors at their global minimizer. However, authors in \cite{RF06a} define the support vectors at the KKT point of $L_r$-SVM. In the following, we define the $L_r$ support vectors  for $L_r$-SVM at the P-stationary point, which is also the local minimizer of $L_r$-SVM.

\vskip 2mm
Summarizing the above analysis, we obtain the following interesting result.

\begin{theorem}{\rm\label{support-vector}  Let ($\bfw^{*};b^{*};\bfu^{*}$) with $\bm{\lambda}^*\in \R^m$ and $\gamma>0$ be a P-stationary  point of (\ref{RR-SVM}) for $\gamma C\geq2$. Then all $L_r$ support vectors must fall into the support hyperplanes $\langle \bfw^*,\bfx\rangle+b^*= \pm1$, which possesses the same feature as the one of hard margin SVM.}
\end{theorem}
\section{Conclusions}
%
In this paper, with the help of explicit expression of proximal operator for ramp loss, we have introduced and characterized a novel first-order necessary and sufficient optimality conditions of $L_r$-SVM. We have defined the $L_r$ support vectors based on the concept of P-stationary point and showed that all $L_r$ support vectors fall into the support hyperplanes for $\gamma C\geq2$. Based on the advance of P-stationarity, could we design an efficient algorithm for $L_r$-SVM? We leave this topic to be investigated in the future.

\vskip 2mm
\noindent{\bf\large Acknowledgements}\\

This work is supported by the National Natural Science Foundation of China (11926348-9, 61866010, 11871183), and the Natural Science Foundation of Hainan Province (118QN181).

\renewcommand{\refname}{{\Large References}}

\end{document}